\documentclass{article}
\usepackage{amsmath}
\usepackage{amsfonts}

\advance \textheight by \topskip
\advance \textheight by 1pt
\makeatother
\def\bea{\begin{eqnarray}}
\def\ena{\end{eqnarray}}
\def\non{\nonumber}
\def\ep{\epsilon}
\def\ot{\otimes}
\def\lar{\longrightarrow}

\def\deg{\hbox{deg}}

\def\Ker{\hbox{Ker}}
\def\tr{\hbox{tr}}

\def\ch{\hbox{ch}}

\def\lar{\longrightarrow}
\newcommand{\qed}{\hbox{\rule{6pt}{6pt}}}
\newtheorem{prop}{Proposition}
\newtheorem{theorem}{Theorem}

\newtheorem{lemma}{Lemma}
\newtheorem{cor}{Corollary}

\newcommand{\qbc}[2]{
\left[
\begin{array}{c}{#1}\\{#2}\end{array}
\right]_q}
\newcommand{\qibc}[2]{
\left[
\begin{array}{c}{#1}\\{#2}\end{array}
\right]_{q^{-1}}}

\title{
The Chiral Space of Local Operators in $SU(2)$-Invariant Thirring Model
}

\author{
Atsushi Nakayashiki\thanks{
Faculty of Mathematics,
Kyushu University,
Ropponmatsu 4-2-1, Fukuoka 810-8560, Japan, \quad 
e-mail: 6vertex@math.kyushu-u.ac.jp}
}

\date{}
\begin{document}
\maketitle
\begin{abstract}
The space of local operators in the $SU(2)$ invariant Thirring model 
($SU(2)$ ITM) is studied by the form factor bootstrap method.
By constructing sets of form factors explicitly we define a susbspace
of operators which has the same character as the level one integrable highest
weight representation of $\widehat{sl_2}$. This makes a correspondence
between this subspace and the chiral space of local operators
 in the underlying conformal field theory, the $su(2)$ Wess-Zumino-Witten model
at level one.
\end{abstract}

\section{Introduction}
The $SU(2)$ invariant Thirring model ($SU(2)$ ITM) 
is a massive integrable model
in two dimensional quantum field theory which can be considered 
as a deformation of the level one $su(2)$ Wess-Zumino-Witten (WZW) model
in conformal field theory. It is expected that there is a one to one
correspondence between local operators in a CFT model and
those in a massive deformation of it \cite{Z}.
In particular there should be a subspace of operators in $SU(2)$ ITM 
which is isomorphic to the level one integrable highest weight 
representation of $\widehat{sl_2}$, the chiral space of the 
level one $su(2)$ WZW model. In this paper we single out such a space
by constructing sets of form factors explicitly.

Let $V=\mathbb{C}v_{+}\oplus \mathbb{C} v_{-}$ be the two dimensional
vector space on which $sl_2$ acts. The $S$-matrix of $SU(2)$ ITM 
is the operator $S(\beta):\,V^{\ot 2}\lar V^{\ot 2}$, the explicit form 
of which is given in (\ref{s-matrix}). 
The set of form factors $(F_{2n})_{n=0}^\infty$, $F_{2n}$ being a $V^{\ot 2n}$-valued
meromorphic function, of a local operator in $SU(2)$ ITM satisfies the
following system of equations:
\bea
&&
P_{i,i+1}S_{i,i+1}(\beta_i-\beta_{i+1})
F_{2n}(\beta_1,\cdots,\beta_{2n})
=F_{2n}(\cdots,\beta_{i+1},\beta_i,\cdots),
\label{iff-1}
\\
&&
P_{2n-1,2n}P_{2n-1,2n-2}\cdots P_{1,2}
F_{2n}(\beta_1-2\pi i,\cdots,\beta_{2n})
=
(-1)^{n}F_{2n}(\beta_2,\cdots,\beta_{2n},\beta_1),
\label{iff-2}
\\
&&
2\pi i \text{res}_{\beta_{2n}=\beta_{2n-1}+\pi i} 
F_{2n}(\beta_1,\cdots,\beta_{2n}) 
\non \\
&& \!\!\!\! 
{}=
\Big(
I-(-1)^{n-1}
S_{2n-1, 2n-2}(\beta_{2n-1}-\beta_{2n-2})\cdots 
S_{2n-1, 1}(\beta_{2n-1}-\beta_{1})
\Big)
\non
\\
&&
\quad
F_{2n-2}(\beta_1,\cdots,\beta_{2n-2})\otimes \mathbf{e}_0,
\label{iff-3}
\ena
where $\mathbf{e}_0=v_{+}\otimes v_{-}-v_{-}\otimes v_{+}$, 
$P:V^{\ot 2}\lar V^{\ot 2}$ is the permutation operator and
$S_{ij}(\beta)$ acts on $i$-th and $j$-th components of $V^{\ot 2n}$
as $S(\beta)$ etc. Conversely any solution $(F_{2n})_{n=0}^\infty$ 
of (\ref{iff-1}), (\ref{iff-2}), (\ref{iff-3}) determines a local operator 
such that its $2n$-particle form factor is $F_{2n}$ \cite{Smir2}. 
Thus the space of local operators can be identified with the space of 
solutions of (\ref{iff-1}), (\ref{iff-2}), (\ref{iff-3}).

The complete description of the solution space of 
(\ref{iff-1}) and (\ref{iff-2}) is known (see \cite{N} and references therein).
Let $R_n$ be the ring of symmetric polynomials
of $x_1,\cdots,x_n$ with $x_j=\exp(\beta_j)$ and $C_n$ the field of 
symmetric, $2\pi i$-periodic and meromorphic functions of 
$\beta_1,\cdots,\beta_n$.
For $\ell\leq n$ consider a polynomial 
$P=P(X_1,\cdots,X_{\ell}|x_1,\cdots,x_{2n})$ such that it is anti-symmetric
in $X_a$'s with the coefficients in $R_{2n}$ satisfying 
$\deg_{X_a}\, P\leq 2n-1$ for any $a$. 
To each $P$ a solution $\Psi_P$ of 
(\ref{iff-1}), (\ref{iff-2}) is constructed through multi-dimensional 
$q$-hypergeometric integrals. It takes the value in the space of 
$sl_2$ highest weight vectors with the weight $2n-2\ell$.
Since $S(\beta)$ and $P$ (the permutation operator) commute with $sl_2$, 
other solutions are created from highest weight solutions
by the actions of $sl_2$ and $C_{2n}$.

It is convenient to describe the space of polynomials as above by some exterior
space.
Let $H^{(2n)}=\oplus_{j=0}^{2n-1}R_{2n}X^j$ and 
$H^{(2n)}_{C_{2n}}=\oplus_{j=0}^{2n-1}C_{2n}X^j$.
Then $\wedge^\ell H^{(2n)}_{C_{2n}}$ can be identified with the space 
of polynomials satisfying the conditions above and 
$\wedge^\ell H^{(2n)}$ becomes a subspace of it.

For each $n\geq 0$ we classify solutions of (\ref{iff-1})-(\ref{iff-3}) 
into the solutions of the form $(0,\cdots,0,F_{2n},F_{2n+2},\cdots)$.
Then Equation (\ref{iff-3}) implies
\bea
&&
\text{res}_{\beta_{2n}=\beta_{2n-1}+\pi i} 
F_{2n}(\beta_1,\cdots,\beta_{2n}) =0.
\label{iff-4}
\ena
If we write $F_{2n}=\Psi_P$, Equation (\ref{iff-4}) should be rewritten as
equations for a polynomial $P$. In fact the following equations imply 
(\ref{iff-4}) 
\bea
&&
P|_{(x_{2n-1},x_{2n})=(x,-x), X_{\ell}=\pm x^{-1}}=0.
\label{nre}
\ena
The converse is conjectured to be true.
Let
\bea
&&
U_{2n,\ell}=\{P\in \wedge^\ell H^{(2n)}|\,\,\text{$P$ is a solution to (\ref{nre})}\}.
\non
\ena
In \cite{N} $U_{2n,\ell}$ is proved to be a free $R_{2n}$-module and its
basis has been constructed explicitly.

To obtain the solution space of (\ref{iff-1}), (\ref{iff-2}), (\ref{iff-4})
one has to consider $U_{2n,\ell}$ modulo polynomials $P$ such that $\Psi_P$
vanishes identically. Such polynomials have been determined completely 
\cite{Tar}.
Namely there are polynomials $\Xi^{(2n)}_1\in U_{2n,1}$ and  
$\Xi^{(2n)}_2\in U_{2n,2}$ such that $\Psi_P$ is identically zero
if and only if $P\in \Xi^{(2n)}_1\wedge^{\ell-1} H^{(2n)}+
\Xi^{(2n)}_2\wedge^{\ell-2} H^{(2n)}$. 
We define
\bea
&&
M_{2n,\ell}=
\frac{U_{2n,\ell}}
{\Xi^{(2n)}_1 \wedge U_{2n,\ell-1}  +\Xi^{(2n)}_2\wedge U_{2n,\ell-2}},
\non
\ena
and, for $\lambda\geq 0$, 
$M^{(0)}_\lambda=\oplus_{2n-2\ell=\lambda}M_{2n,\ell}$.
Then $M^{(0)}_\lambda$ becomes isomorphic to the space of 
$sl_2$ highest weight vectors in the level one integrable highest 
weight representation $V(\Lambda_0)$ of $\widehat{sl_2}$ \cite{N}.

In this paper, to each $P_{2m}\in U_{2m,\ell}$, we construct polynomials
$P_{2n}$, $n>m$ such that $(\Psi_{P_{2n}})_{n=0}^\infty$ satisfies 
(\ref{iff-3}) (Theorem \ref{main}), where we set $P_{2n}=0$ for $n<m$ . 
This set of functions specifies a local operator which is a highest weight 
vector of $sl_2$.
The operators corresponding to non-highest weight vectors are obtained 
from the ones above by the action of $sl_2$.
Let $V_\lambda$ be the finite dimensional irreducible representation
of $sl_2$ with the highest weight $\lambda$.
Then $\oplus_\lambda M^{(0)}_\lambda \otimes V_\lambda$ specifies
the subspace of local operators in $SU(2)$ ITM which is isomorphic to 
$V(\Lambda_0)$.

Finally we briefly remark on the previous results on
the determination of the space of local operators in massive integrable models.
For several models with diagonal S-matrices the characters of the space of 
initial form factors, which are the first non-vanishing form factors, have
been calculated \cite{K}. They are shown to take fermionic forms \cite{KMM} 
of the corresponding CFT characters. 
For Sine-Gordon model and SU(2) ITM, which are two of simple models with
non-diagonal S-matrices, important results toward determining all local 
opertators have been obtained in \cite{Smir3} and \cite{Luk}.
However the construction of $2n$-particle form factors for every $n$ 
to each element of the underlying CFT characters has never been done even
for models with diagonal S-matrices other than Ising model \cite{CM,Ch}
and the model considered in \cite{BBS}.
We solve this problem for SU(2) ITM under certain assumptions.

The plan of the present paper is as follows.
In section 2 we review the results of \cite{N}, in particular, 
the description of the basis of $U_{2n,\ell}$.
The formulae for polynomials $P_{2n}$, $n\geq m$ and 
consequently for form factors are given in section 3.
In section 4 some theorems on symmetric polynomials which is
used in section 3 are proved. They are interesting by itself.
The form factors of local operators in the anti-chiral subspace
are given in section 5.
 In appendix A the definitions and
properties of the function $\zeta(\beta)$ which appeared in the formulae of
form factors is briefly explained for the sake of readers' convenience.

\section{Space of minimal form factors}
In this section we recall the results of \cite{N}.

Let $V=\mathbb{C}v_{+}\oplus \mathbb{C} v_{-}$ be the vector representation
of $sl_2=\mathbb{C}e\oplus\mathbb{C}h\oplus\mathbb{C}f$ with the correspondence
$v_{+}={}^t(1,0)$, $v_{-}={}^t(0,1)$ and $P: V^{\otimes2}\lar V^{\otimes2}$
the permutation operator $P(v_{\ep_1}\ot v_{\ep_2})=v_{\ep_2}\ot v_{\ep_1}$.
The $S$-matrix of the $SU(2)$ invariant Thirring model is the linear operator
in $End(V^{\otimes2})$ given as
\bea
&&
S(\beta)=S_0(\beta)\hat{S}(\beta),
\quad
\hat{S}(\beta)=
\frac{\beta-\pi i P}{\beta-\pi i},
\quad
S_0(\beta)=
\frac{\Gamma(\frac{\pi i+\beta}{2\pi i})\Gamma(\frac{-\beta}{2\pi i})}
{\Gamma(\frac{\pi i-\beta}{2\pi i})\Gamma(\frac{\beta}{2\pi i})}.
\label{s-matrix}
\ena
For a linear operator $A$ in $End(V^{\otimes 2})$ we define 
$A_{ij}\in End(V^{\otimes n})$ as follows.
Write
\bea
&&
A=\sum_a B_a\otimes C_a,
\quad
B_a,C_a\in \text{End}(V).
\non
\ena
Then
\bea
&&
A_{ij}=\sum_a\, 
1\otimes\cdots\otimes B_a\otimes\cdots \otimes C_a\otimes \cdots \otimes1,
\non
\ena
where
$B_a$ and $C_a$ are in i-th and j-th components respectively.

Consider the following system of equations for a $V^{\otimes n}$-valued 
function $F$:
\bea
&&
P_{i,i+1}S_{i,i+1}(\beta_i-\beta_{i+1})
F(\beta_1,\cdots,\beta_n)
=F(\cdots,\beta_{i+1},\beta_i,\cdots),
\label{feq-1}
\\
&&
P_{n-1,n}P_{n-1,n-2}\cdots P_{1,2}
F(\beta_1-2\pi i,\cdots,\beta_n)
=
(-1)^{\frac{n}{2}}F(\beta_2,\cdots,\beta_n,\beta_1).
\label{feq-2}
\ena
All meromorphic solutions of those equations are given by the multi-dimensional
integrals of $q$-hypergeometric type in the following manner.

Notice that $S(\beta)$ commutes with the action of $sl_2$. Therefore any 
solution can be obtained, by the action of $sl_2$, from the solutions
taking values in the space of highest weight vectors 
$(V^{\ot n})^{sing}_\lambda$:
\bea
&&
(V^{\ot n})^{sing}_\lambda=\{v\in V^{\ot n}\,\vert\, ev=0,\, hv=\lambda v\}.
\non
\ena
We remark that $(V^{\ot n})^{sing}_\lambda\neq \{0\}$ if and only if
$\lambda$ is written as $\lambda=n-2\ell$ for some $0\leq \ell \leq n/2$.

For a subset $M\subset \{1,2,...,n\}$ we set
\bea
&&
v_M=v_{\epsilon_1}\otimes\cdots\otimes v_{\epsilon_n}\in V^{\ot n},
\quad
M=\{j\,|\, \epsilon_j=-\,\}.
\non
\ena

Let $C_n$ be the field of symmetric, $2\pi i$-periodic and meromorphic 
functions of $\beta_1,...,\beta_n$ and 
\bea
&&
H^{(n)}:=\oplus_{k=0}^{n-1} R_{n} X^k,
\quad
H^{(n)}_{C_n}:=\oplus_{k=0}^{n-1} C_{n} X^k
\non
\ena
the space of polynomials of $X$ of degree at most $n-1$ with the coefficients
in $R_n$ and $C_n$ respectively. We identify the $\ell$-th exterior product 
space
$\wedge^\ell H^{(n)}_{C_n}$ with the space of anti-symmetric polynomials of
$X_1,...,X_\ell$ with the coefficients in $C_n$ by
\bea
X^{i_1}\wedge\cdots\wedge X^{i_\ell}
\mapsto
\text{Asym}(X_1^{i_1}\cdots X_\ell^{i_\ell})
=\sum_{\sigma\in S_\ell}\text{sgn}\,\sigma\, X_{\sigma(1)}^{i_1}\cdots 
X_{\sigma(\ell)}^{i_\ell}.
\non
\ena
We use the variables $X_a$, $x_j$ and $\alpha_a$, $\beta_j$ which are related
by $X_a=\exp(-\alpha_a)$, $x_j=\exp(\beta_j)$.

Now for each $P\in \wedge^\ell H^{(n)}_{C_n}$ define the $V^{\otimes n}$-valued function $\psi_P$ by
\bea
&&
\psi_P=\sum_{\sharp M=\ell} I_M(P)v_M,
\non
\ena
where
\bea
&&
I_M(P):=
\int_{C^\ell}\prod_{a=1}^\ell d\alpha_a \prod_{a=1}^\ell\phi_n(\alpha_a)
w_M
\frac{P(X_1,\cdots,X_\ell|x_1,\cdots,x_n)}
{\prod_{a=1}^\ell\prod_{j=1}^{n}(1-X_ax_j)},
\label{eq2-1}
\ena
where, $w_M=Asym(g_M)$ and, for 
$M=(m_1,\cdots,m_\ell)$ with $m_1<\cdots<m_\ell$,
\bea
&&
\phi_n(\alpha)=\prod_{j=1}^n
\frac{\Gamma(\frac{\alpha-\beta_j+\pi i}{-2\pi i})}
{\Gamma(\frac{\alpha-\beta_j}{-2\pi i})},
\quad
g_M=\prod_{a=1}^\ell
\Big(
\frac{1}{\alpha_a-\beta_{m_a}}
\prod_{j=1}^{m_a-1}
\frac{\alpha_a-\beta_j+\pi i}{\alpha_a-\beta_j}
\Big)
\prod_{1\leq a<b\leq \ell}(\alpha_a-\alpha_b+\pi i).
\non
\ena
The integration contour $C$ is defined as follows.
It goes from $-\infty$ to $+\infty$, separating two sets 
$\{\beta_j-2\pi i k\,|\,1\leq j\leq n,\, k\geq 0\}$ and
$\{\beta_j+(2k-1)\pi i \,|1\leq j\leq n,\, k\geq 0\}$.

We set
\bea
&&
\Psi_P=\exp\Big(\frac{n}{4}\sum_{j=1}^n\beta_j\Big)
\prod_{j<k}\zeta(\beta_j-\beta_k)\psi_P,
\ena
where $\zeta(\beta)$ is a certain meromorphic function described 
in Appendix A.

Let $Sol^n_\lambda$ be the 
space of meromorphic solutions of (\ref{feq-1}) and (\ref{feq-2})
taking values in $(V^{\ot n})^{sing}_{\lambda}$.
Then $\Psi_P\in Sol^n_{n-2\ell}$. 

\vskip2mm
\noindent
{\bf Remark.} The integral $I_M(P)$ is well defined for any polynomial
$P(X_1,...,X_\ell|x_1,...,x_n)$ such that $\deg_{X_a}\, P\leq n$ for any
$a$. If $P$ is
symmetric in $x_j$'s, the function $\Psi_P$ becomes a solution of
(\ref{feq-1}) and (\ref{feq-2}). Due to the anti-symmetry of $w_M$
the relation $\ell!\Psi_P=\Psi_{Asym(P)}$ holds. The highest degree $X_a^n$
can be eliminated by $\text{Ker}\, \Psi$ below. Thus for any $P$ we have 
$\Psi_P=\Psi_{P'}$ for some $P'\in \wedge^\ell H^{(n)}_{C_n}$.
We later use this property in the description of polynomials corresponding to
form factors.
\vskip2mm

We have the linear map:
\bea
&&
\Psi: \wedge^\ell H^{(n)}_{C_n}\lar Sol^n_{n-2\ell},
\quad
P\mapsto \Psi_P.
\non
\ena
This map is surjective.
The kernel of $\Psi$ is described as follows.
Let $\Theta^{(n)}_{\pm}(X)=\prod_{j=1}^n(1\pm Xx_j)$ and
\bea
&&
2\Xi_1^{(n)}(X)=\Theta^{(n)}_{+}(X)+(-1)^{n-1}\Theta^{(n)}_{-}(X),
\non
\\
&&
2\Xi_2^{(n)}(X_1,X_2)=
\Big(\Theta^{(n)}_{+}(X_1)\Theta^{(n)}_{+}(X_2)
-
\Theta^{(n)}_{-}(X_1)\Theta^{(n)}_{-}(X_2)
\Big)
\frac{X_1-X_2}{X_1+X_2}
\non
\\
&&
\qquad\qquad\qquad
+(-1)^n
\Big(
\Theta^{(n)}_{+}(X_1)\Theta^{(n)}_{-}(X_2)
-
\Theta^{(n)}_{+}(X_2)\Theta^{(n)}_{-}(X_1)
\Big).
\non
\ena
Then
\bea
&&
\Ker\, \Psi=\wedge^{\ell-1} H^{(n)}_{C_n}\wedge \Xi_1^{(n)}
+
\wedge^{\ell-2} H^{(n)}_{C_n}\wedge \Xi_2^{(n)}.
\non
\ena

Consider the restriction of $\Psi$ to $\wedge^\ell H^{(n)}$. 
Then the function $\Psi_P$ has at most a simple pole at 
$\beta_n=\beta_{n-1}+\pi i$ as a function of $\beta_n$.
In \cite{N} the subspace of $\wedge^\ell H^{(n)}$ such that
the function $\Psi_P$ is holomorphic at $\beta_n=\beta_{n-1}+\pi i$ 
has been studied. Let us recall the results.

For a polynomial $P(X_1,...,X_\ell|x_1,\cdots,x_n)$ in $X_a$'s with the coefficients in Laurent polynomials of $x_j$'s define
\bea
&&
\rho_{\pm}(P)=P(X_1,\cdots,X_{\ell-1},\pm x^{-1}|x_1,\cdots,x_{n-2},x,-x),
\non
\ena
and set
\bea
&&
U_{n,\ell}=\{\, P\in \wedge^\ell H^{(n)}\vert \, 
\rho_{+}(P)=\rho_{-}(P)=0\,\}.
\non
\ena
Then $\Xi^{(n)}_1\in U_{n,1}$, $\Xi^{(n)}_2\in U_{n,2}$ and if 
$P$ is in $U_{n,\ell}$ then $\Psi_P$ becomes regular at 
$\beta_n=\beta_{n-1}+\pi i$.
Let 
\bea
&&
M_{n,\ell}=
\frac{U_{n,\ell}}
{U_{n,\ell-1} \wedge \Xi^{(n)}_1+U_{n,\ell-2}\wedge \Xi^{(n)}_2},
\non
\ena
and, for $i=0,1$, $\lambda\in \mathbb{Z}_{\geq 0}$,
\bea
&&
M^{(i)}_{\lambda}=
\oplus_{n\equiv i\, mod.\, 2,\, n-2\ell=\lambda} M_{n,\ell}.
\non
\ena

\vskip2mm
\noindent
{\bf Remark}.
We conjecture the following equation:
\bea
&&
U_{n,\ell}\cap
\Big(
\wedge^{\ell-1} H^{(n)}_{C_n}\wedge \Xi_1^{(n)}
+
\wedge^{\ell-2} H^{(n)}_{C_n}\wedge \Xi_2^{(n)}
\Big)
=
U_{n,\ell-1} \wedge \Xi^{(n)}_1+U_{n,\ell-2}\wedge \Xi^{(n)}_2.
\label{ker-conj}
\ena
For $n=1,2$ (\ref{ker-conj}) holds.
If (\ref{ker-conj}) holds, then the space $M_{n,\ell}$ becomes 
isomorphic to the subspace of $Sol^n_{n-2\ell}$.
\vskip3mm

We define a degree of an element in $U_{n,\ell}$, 
denoted by $\deg_1$, assigning 
$\deg_1\, X_{a}=-1$, $\deg_1\, x_j=1$. 
Then we set $\deg_2\, P=\deg_1\, P+n^2/4$ for a homogeneous element $P$ 
in $U_{n,\ell}$. 
Then $\deg_1\, \Xi^{(n)}_1=\deg_1\, \Xi^{(n)}_2=0$.
We introduce a grading on $M_{n,\ell}$ by $\deg_2$.
In general for a graded vector space $V=\oplus_n V_n$ we define
its character by
\bea
&&
\ch\, V=\sum_n q^n \dim\, V_n.
\non
\ena

Let $\Lambda_i$ be 
the fundamental weights of $\widehat{sl_2}$, $V(\Lambda_i)$ the 
integrable highest weight representation of $\widehat{sl_2}$ with 
the highest weight $\Lambda_i$, 
$V(\Lambda_i|\lambda)$, $\lambda\in \mathbb{Z}_{\geq 0}$, the subspace of 
$V(\Lambda_i)$ consisting of $sl_2$-highest weight vectors with the weight 
$\lambda$ and $d$ the scaling element of $\widehat{sl_2}$ \cite{Kac}.

\begin{theorem}$(\cite{N})$\label{th-n}
\bea
&&
\ch\, M^{(i)}_{\lambda}=
\tr_{V(\Lambda_i|\lambda)}(q^{-d+i/4})=
\sum_{n-2\ell=\lambda,\, n\equiv i\, mod.\, 2}
\frac{q^{\frac{n^2}{4}}}{[n]_q!}
\Big(
\qbc{n}{\ell}-\qbc{n}{\ell-1}
\Big),
\label{branching}
\ena
where
\bea
&&
[n]_q=1-q^n,
\quad
[n]_q!=\prod_{j=1}^n[j]_q,
\quad
\qbc{n}{\ell}=\frac{[n]_q!}{[\ell]_q![n-\ell]_q!}.
\non
\ena
\end{theorem}

The $R_n$-module $U_{n,\ell}$ becomes a free module and a basis is 
given explicitly. Since we consider the case of $n$ even in this paper,
we recall the basis only in the case of $n$ even here.

Let $e^{(n)}_k$ be the elementary symmetric polynomial defined by
\bea
&&
\prod_{j=1}^n(1+x_jt)=\sum_{k=0}^ne^{(n)}_kt^k.
\non
\ena
The symmetric polynomial $P^{(2n)}_{r,s}$, $1\leq r\leq n$, $s\in \mathbb{Z}$ 
is defined by the following recursion relations:
\bea
P^{(2n)}_{1,s}&=&e^{(2n)}_{2s-1},
\non
\\
P^{(2n)}_{r,s}&=&P^{(2n)}_{r-1,s+1}-e^{(2n)}_{2s}P^{(2n)}_{r-1,1}
\quad \hbox{for $r\geq 2$}.
\non
\ena
For a function $f(x_1,\cdots,x_{n})$ of $x_j$'s we denote $\bar{f}$ 
the function obtained from $f$ by specializing the last two variables to
$x,-x$:
\bea
&&
\bar{f}=f(x_1,\cdots,x_{n-2},x,-x), 
\quad
x=x_{n-1},
\non
\ena
if it has a sense.
Then $P^{(2n)}_{r,s}$'s satisfy
\bea
&&
\overline{P^{(2n)}_{r,s}}=P^{(2n-2)}_{r,s}-x^2P^{(2n-2)}_{r,s-1}.
\label{bar-P}
\ena

We set
\bea
&&
v^{(2n)}_0=\sum_{j=0}^ne^{(2n)}_jX^{2j},
\quad
v^{(2n)}_r=\sum_{s=1}^nP^{(2n)}_{r,s}X^{2(s-1)},
\qquad
w^{(2n)}_r=Xv^{(2n)}_r, 
\quad
1\leq r\leq n,
\ena
and, for $1\leq j\leq n$,
\bea
2\xi_j^{(2n)}(X_1,X_2)&=&
\frac{X_1-X_2}{X_1+X_2}
\Big(
v_0^{(2n)}(X_1)w^{(2n)}_j(X_2)+v_0^{(2n)}(X_2)w^{(2n)}_j(X_1)
\Big)
\non
\\
&&
-v_0^{(2n)}(X_1)w^{(2n)}_j(X_2)+v_0^{(2n)}(X_2)w^{(2n)}_j(X_1).
\ena
Equation (\ref{bar-P}) implies that $v^{(2n)}_r,w^{(2n)}_r\in U_{2n,1}$, 
$1\leq r\leq n$ and $\xi_j^{(2n)}\in U_{2n,2}$.

We use the multi-index notations like 
$v_I=v_{i_1}\wedge\cdots\wedge v_{i_{\ell}}$ for $I=(i_1,\cdots,i_\ell)$.
Then, in general, we have the following theorem.

\begin{theorem}\label{basis}
As an $R_{2n}$-module, $U_{2n,\ell}$ is a free module with the following
set of elements as a basis:
\bea
&&
v^{(2n)}_I\wedge w^{(2n)}_J\wedge \xi^{(2n)}_K,
\non
\\
&&
I=(i_1,\cdots,i_{\ell_1}),
\quad
1\leq i_1<\cdots<i_{\ell_1}\leq n,
\non
\\
&&
J=(j_1,\cdots,j_{\ell_2}),
\quad
1\leq j_1<\cdots<j_{\ell_2}\leq n-\ell_1-\ell_3,
\non
\\
&&
K=(k_1,\cdots,k_{\ell_3}),
\quad
1\leq k_1\leq \cdots \leq k_{\ell_3}\leq n-\ell_1-\ell_3+1,
\non
\\
&&
\ell_1+\ell_2+2\ell_3=\ell.
\non
\ena
\end{theorem}

\section{Chiral subspace}
Let $m$ and $r$ be non-negative integers such that $0\leq r\leq m$.
We fix them once for all in this section. We set $\lambda=2m-2r$ and
$\ell_{2n}=r+n-m$. Define $P_{2n}=0$ for $n<m$.
To each $P_{2m}\in U_{2m,r}$ we shall construct polynomials
$P_{2n}(X_1,...,X_{\ell_{2n}}|x_1,...,x_{2n})$, $n\geq m$ such that 
$(\Psi_{P_{2n}})_{n=0}^\infty$ satisfy 
(\ref{iff-1}), (\ref{iff-2}),(\ref{iff-3}).

Let us first discuss the symmetry of equations (\ref{iff-1}), (\ref{iff-2}),
(\ref{iff-3}). Let $p_i=\sum_{j=1}^\infty x_j^i$ and 
$p^{(n)}_i=\sum_{j=1}^nx_j^i$ for $i\neq 0$. Suppose that 
$(F_{2n})_{n=0}^\infty$ is a solution of (\ref{iff-1}), (\ref{iff-2}),
(\ref{iff-3}). Then, for any $s\in \mathbb{Z}$, $(p^{(2n)}_{2s-1}F_{2n})_{n=0}^\infty$ is again a solution of the same equations, 
since $\overline{p^{(2n)}_{2s-1}}=p^{(2n-2)}_{2s-1}$. 
Thus the space of solutions of (\ref{iff-1}), (\ref{iff-2}), (\ref{iff-3})
becomes a module over the ring $\mathbb{Z}[p_{\pm1},p_{\pm3},\cdots]$.
The action of $p_{2s-1}$ is that of the local integral of motion with 
spin $2s-1$. 
We introduce the generating function of $p_{2s-1}^{(n)}$:
\bea
&&
E_{odd}^{(n)}(t)=\exp(\sum_{j\in\mathbb{Z}}t_{2j-1}p^{(n)}_{2j-1}),
\non
\ena
which satisfies the equation 
\bea
&&
\overline{E_{odd}^{(n)}(t)}=E_{odd}^{(n-2)}(t).
\non
\ena

Next we consider the functions
\bea
&&
Q_{2n}^{(\pm)}(z)=\frac{\prod_{a=r+1}^{\ell_{2n}}(1-X_a^{\mp 2}z^2)}
{\prod_{j=1}^{2n}(1-x_j^{\pm1}z)},
\non
\ena
for $n\geq m$.
For $n=m$ the empty product in the numerator is understood as $1$, that is,
\bea
&&
Q_{2m}^{(\pm)}(z)=\frac{1}{\prod_{j=1}^{2m}(1-x_j^{\pm1}z)}=
\sum_{j=0}^\infty h^{(2m)\pm}_jz^j,
\non
\ena
where $h^{(2m)}_j:=h^{(2m)+}_j$ is the complete symmetric function \cite{Mac}.
They satisfy the equation
\bea
&&
\rho_{\pm}\Big(Q_{2n}^{(\ep)}(z)\Big)=
Q_{2n-2}^{(\ep)}(z),
\label{bar-Q}
\ena
for $\ep=\pm$.
These functions were extracted from the integral formulae for the form factors
of local operators defined by Lukyanov \cite{NT,NPT,Luk}.

Take $I,J,K$ such that they satisfy conditions in Theorem \ref{basis} with
$\ell_1+\ell_2+2\ell_3=r$. We set
\vskip2mm
\bea
&&
P_{2m}=E_{odd}^{(2m)}(t)\prod_{i=1}^{m}Q_{2m}^{(+)}(z_i)
v^{(2m)}_I\wedge w^{(2m)}_J\wedge \xi^{(2m)}_K,
\label{ff-1}
\\
&&
P_{2n}=c_{2n}E_{odd}^{(2n)}(t)\prod_{i=1}^{m}Q_{2n}^{(+)}(z_i)
v^{(2n)}_I\wedge w^{(2n)}_J\wedge \xi^{(2n)}_K
\prod_{a=r+1}^{\ell_{2n}}X_a^{2n+1+2r-2a},
\quad
\text{for $n>m$},
\label{ff-2}
\ena
\vskip2mm
\noindent
where the constant $c_{2n}$ is given by
\bea
&&
c_{2n}=(-1)^{\frac{1}{2}(n-m)(n+m-2r-1)}
(i\zeta(-\pi i))^{m-n}
(2\pi i)^{(m-n)(n+r)}.
\label{const-c}
\ena
In (\ref{ff-1}), (\ref{ff-2}), $v^{(2n)}_I\wedge w^{(2n)}_J\wedge \xi^{(2n)}_K$
is understood as the polynomial of $X_1,\cdots,X_r$.
We expand $P_{2n}$ as
\bea
&&
P_{2n}=\sum_{\alpha,\gamma}P_{2n,\alpha,\gamma}t^\alpha z^\gamma,
\non
\ena
where $\alpha=(...,\alpha_{-1},\alpha_1,\alpha_3,...)$, 
$\gamma=(\gamma_1,...,\gamma_m)$ and $t_\alpha=\prod t_i^{\alpha_i}$,
$z\gamma=\prod z_i^{\gamma_i}$.

\begin{theorem}\label{main}
$(1)$. The set of functions $(\Psi_{P_{2n}})_{n=0}^\infty$ satisfies
$(\ref{iff-1})$, $(\ref{iff-2})$, $(\ref{iff-3})$ and each $\Psi_{P_{2n}}$
takes the value in $(V^{\ot 2n})^{sing}_{\lambda}$, $\lambda=2m-2r$.

\noindent
$(2)$. The following property holds for all $n$:
\bea
&&
\Psi_{P_{2n,\alpha,\gamma}}(\beta_1+\theta,\cdots,\beta_{2n}+\theta)
=\exp(\theta\deg_2\, P_{2m,\alpha,\gamma})
\Psi_{P_{2n,\alpha,\gamma}}(\beta_1,\cdots,\beta_{2n}).
\non
\ena
\end{theorem}

The property (2) in Theorem \ref{main} was proved in \cite{N,NT}. 
It means that the Lorentz spin of the local operator
specified by the set of form factors $(\Psi_{P_{2n}})_{n=0}^\infty$
is $\deg_2\, P_{2m,\alpha,\gamma}$. 
Explicitly $\deg_2\, P_{2m,\alpha,\gamma}$ is given by
\bea
&&
\deg_2\, P_{2m,\alpha,\gamma}=m^2+\sum i\,\alpha_i+\sum \gamma_i+
\deg_1(v^{(2m)}_I\wedge w^{(2m)}_J\wedge \xi^{(2m)}_K).
\non
\ena

As shown in Theorem \ref{th-symfunc} in the next section
$R_{2m}$ is a free $\mathbb{C}[p^{(2m)}_1,p^{(2m)}_3,\cdots,p^{(2m)}_{2m-1}]$
module with the $n$-th ordered products of $h_{2j}$'s as a basis, that is, 
\bea
&&
R_{2m}=\oplus_{0\leq r_1\leq \cdots\leq r_m}
\mathbb{C}[p^{(2m)}_1,p^{(2m)}_3,\cdots,p^{(2m)}_{2m-1}]
h^{(2m)}_{2r_1}\cdots h^{(2m)}_{2r_m}.
\non
\ena
Thus any element of $U_{2m,r}$ can be expressed as a linear combination
of $P_{2m,\alpha,\gamma}$'s, that is,
\bea
&&
U_{2m,r}=\sum_{\alpha,\gamma}\mathbb{C}P_{2m,\alpha,\gamma},
\non
\ena
where the summation is taken for all $\gamma$ and $\alpha$ 
with $\alpha_i=0$, $i<0$ or $i>2m$. 

Let ${\cal O}^{sing}_\lambda(2m)$ be the space of meromorphic 
solutions to (\ref{iff-1}),
(\ref{iff-2}), (\ref{iff-3}) such that $F_{2n}=0$, $n<m$ and 
$F_{2n}\in Sol^{2n}_\lambda$ for all $n$. 
We set ${\cal O}^{sing}_\lambda={\cal O}^{sing}_\lambda(0)$. 
Then we have
\bea
&&
{\cal O}^{sing}_\lambda={\cal O}^{sing}_\lambda(0)
\supset
{\cal O}^{sing}_\lambda(2)
\supset
{\cal O}^{sing}_\lambda(4)
\supset\cdots.
\non
\ena
For each $m$ there map is a map
\bea
&&
\Psi: U_{2m,r}\lar {\cal O}^{sing}_{2m-2r}(2m),
\quad
P_{2m,\alpha,\gamma}\mapsto (\Psi_{P_{2n,\alpha,\gamma}})_{n=0}^\infty,
\non
\ena
where $\alpha$ is such that $\alpha_i=0$ for $i<0$ or $i>2m$.
It induces the map
\bea
&&
\Psi: M_{2m,r}\lar {\cal O}^{sing}_{2m-2r}(2m).
\non
\ena
Finally we obtain the map
\bea
&&
\Psi:M^{(0)}_\lambda=\oplus_{\lambda=2m-2r} M_{2m,r}\lar
{\cal O}^{sing}_{\lambda}.
\label{schiral-map}
\ena

Let ${\cal O}$ be the space of meromorphic solutions of 
(\ref{iff-1}), (\ref{iff-2}), (\ref{iff-3}) and 
$V_\lambda\simeq \mathbb{C}^{\lambda+1}$ the
irreducible resresentation of $sl_2$ with the highest weight $\lambda$. 
We consider $M^{(0)}_\lambda$ as the trivial $sl_2$-module. 
Then (\ref{schiral-map}) induces the map of $sl_2$-modules
\bea
&&
\oplus_\lambda M^{(0)}_\lambda\otimes V_{\lambda}
\lar
{\cal O}.
\label{chiral-map}
\ena
We call the image of (\ref{chiral-map}) the chiral subspace of the space of 
local operators in SU(2) ITM. 
If we assume the conjecture (\ref{ker-conj}), the map (\ref{chiral-map}) 
is injective and the chiral subspace is isomorphic to $V(\Lambda_0)$ 
due to Theorem \ref{th-n}.

Now let us prove Theorem \ref{main}.
The following theorem has been proved in \cite{NT}.

\begin{theorem}$(\cite{NT})$\label{nt-th}
Let $P_{2n}(X_1,...,X_{\ell_{2n}}|x_1,...,x_{2n})$, $n\geq m$ be polynomials 
in $X_1$, ..., $X_{\ell_{2n}}$ with the coefficients in the ring of 
symmetric Laurent polynomials of $x_j$'s such that 
$\deg_{X_a}\, P\leq 2n$
for $1\leq a\leq \ell_{2n}$. Set $P_{2n}=0$ for $n<m$. 
Suppose that $(P_{2n})_{n=0}^\infty$ satisfy the following
conditions: there exists a set of polynomials $(P_{2n}')$ such that
\bea
&&
Asym\Big(\overline{P_{2n}}\Big)=
Asym\Big(\prod_{a=1}^{\ell_{2n}-1}(1-x^2X_a^2)P_{2n}'\Big),
\label{cond-1}
\\
&&
P_{2n}'|_{X_{\ell_{2n}}=\pm x^{-1}}=\pm x^{-(2n-1)}d_{2n}P_{2n-2},
\label{cond-2}
\\
&&
d_{2n}=\frac{2\pi}{\zeta(-\pi i)}(-2\pi i)^{-r+m-2n},
\ena
where $Asym$ is the anti-symmetrization with respect to 
$X_{r+1},...,X_{\ell_{2n}}$. 
Then $(\Psi_{P_{2n}})_{n=0}^\infty$ satisfies 
$(\ref{iff-1})$, $(\ref{iff-2})$, $(\ref{iff-3})$.
\end{theorem}

Let us show that (\ref{ff-1}), (\ref{ff-2}) satisfy
(\ref{cond-1}) and (\ref{cond-2}).

\begin{lemma}\label{rec-v}
For $0\leq i\leq n$, $1\leq j\leq n$ we have
\bea
&&
\overline{v^{(2n)}_i(X)}=(1-x^2X^2)v^{(2n-2)}_i(X),
\quad
\overline{w^{(2n)}_j(X)}=(1-x^2X^2)w^{(2n-2)}_j(X),
\non
\\
&&
\overline{\xi^{(2n)}_j(X_1,X_2)}=
\prod_{a=1}^2(1-x^2X^2_a)\xi^{(2n-2)}_j(X_1,X_2).
\non
\ena
\end{lemma}

\vskip2mm
\noindent
{\it Proof.}
Using (\ref{bar-P}) we easily have the results for $v^{(2n)}_i$ 
and $w^{(2n)}_i$ for $i\neq 0$. The equation for $v^{(2n)}_0$
is similarly verified using 
\bea
&&
\overline{e^{(2n)}_k}=e^{(2n-2)}_k-x^2e^{(2n-2)}_{k-2}.
\non
\ena
The equation for $\xi^{(2n)}_j$ follows from those of
$v^{(2n)}_i$ and $w^{(2n)}_i$. \qed

By this lemma we have 
\bea
&&
\overline{v^{(2n)}_I\wedge w^{(2n)}_J\wedge \xi^{(2n)}_K}=
\prod_{a=1}^r(1-x^2X_a^2)\,
v^{(2n-2)}_I\wedge w^{(2n-2)}_J\wedge \xi^{(2n-2)}_K.
\label{bar-vwxi}
\ena

The following lemma is proved in \cite{NT}.

\begin{lemma}\label{lemma-nt}
Let $Asym$ denote the anti-symmetrization with respect to 
$X_{r+1},...,X_{\ell_{2n}}$. We have
\bea
&&
Asym\big(
\prod_{a=r+1}^{\ell_{2n}}X_a^{2n+1+2r-2a}
\big)
\non
\\
&&
=
(-1)^{\ell_{2n}-r-1}
Asym\big(
X_{\ell_{2n}}^{2n-1}
\prod_{a=r+1}^{\ell_{2n}-1}(1-x^2X_a^2)X_a^{2n-1+2r-2a}
\big).
\label{asym-1}
\ena
\end{lemma}

It follows from (\ref{bar-vwxi}) that
\bea
&&
\overline{P_{2n}}=\prod_{a=1}^r(1-x^2X_a^2)P_{2n}'',
\non
\\
&&
P_{2n}''=c_{2n}E_{odd}^{(2n-2)}(t)\prod_{i=1}^m\overline{Q_{2n}^{(+)}(z_i)}
v^{(2n-2)}_I\wedge w^{(2n-2)}_J\wedge \xi^{(2n-2)}_K
\prod_{a=r+1}^{\ell_{2n}}X_a^{2n+1+2r-2a}.
\non
\ena
Notice that $\prod_{i=1}^mQ_{2n}(z_i)$ is symmetric 
with respect to $X_{r+1},...,X_{\ell_{2n}}$ and 
$v^{(2n-2)}_I\wedge w^{(2n-2)}_J\wedge \xi^{(2n-2)}_K$ does not contain
$X_{r+1},...,X_{\ell_{2n}}$. Then using Lemma \ref{lemma-nt} we have
\bea
Asym(\overline{P_{2n}})&=&
\prod_{a=1}^r(1-x^2X_a^2)Asym(P_{2n}'')
\non
\\
&=&
Asym\big(\prod_{a=1}^{\ell_{2n}-1}(1-x^2X_a^2)P_{2n}'\big),
\non
\ena
where
\bea
P_{2n}'&=&
(-1)^{\ell_{2n}-r-1}
c_{2n}
E_{odd}^{(2n-2)}(t)\prod_{i=1}^m\overline{Q_{2n}^{(+)}(z_i)}
v^{(2n-2)}_I\wedge w^{(2n-2)}_J\wedge \xi^{(2n-2)}_K
\non
\\
&&
\times
X_{\ell_{2n}}^{2n-1}\prod_{a=r+1}^{\ell_{2n}-1}X_a^{2n-1+2r-2a}.
\non
\ena
Due to (\ref{bar-Q}) we have
\bea
P_{2n}'|_{X_{\ell_{2n}}=\pm x^{-1}}&=&\pm
(-1)^{n-m-1}
c_{2n}
E_{odd}^{(2n-2)}(t)\prod_{i=1}^mQ_{2n-2}^{(+)}(z_i)
v^{(2n-2)}_I\wedge w^{(2n-2)}_J\wedge \xi^{(2n-2)}_K
\non
\\
&&\times
x^{-(2n-1)}\prod_{a=r+1}^{\ell_{2(n-1)}}X_a^{2n-1+2r-2a}
\non
\\
&=&
\pm(-1)^{n-m-1}\frac{c_{2n}}{c_{2n-1}}x^{-(2n-1)}P_{2n-2}.
\non
\ena
Finally the recursion relation
\bea
&&
\frac{c_{2n}}{c_{2n-1}}=(-1)^{n-m-1}d_{2n},
\quad
c_{2m}=1,
\non
\ena
uniquely determines $c_{2n}$, $n\geq m$ and the explicit form of which is given
by (\ref{const-c}). \qed

\section{Some properties of symmetric functions}
In this section we omit $(n)$ of $e_j^{(n)}$, $p_j^{(n)}$, $h_j^{(n)}$
if there is no fear of confusion on the number of the variables.

Let $R_n^{odd}=
\mathbb{C}[p_{2j+1}|0\leq j\leq (n-1)/2]\subset R_n$ 
and $n'=[n/2]$ be the largest integer not exceeding $n/2$.

\begin{theorem}\label{th-symfunc}
The module $R_n$ is a free $R_n^{odd}$-module with the elements
$\{h_{2i_1}\cdots h_{2i_{n'}}|0\leq i_1\leq\cdots\leq i_{n'}\}$ as a basis:
\bea
&&
R_n=\oplus_{0\leq i_1\leq\cdots\leq i_{n'}}
R_n^{odd}h_{2i_1}\cdots h_{2i_{n'}}.
\non
\ena
\end{theorem}

For the proof of this theorem we need to prepare some propositions and lemmas.
We denote $R_n'$ the subring of $R_n$ generated by 
$e_j$ with $j$ odd, $j\leq n$:
\bea
&&
R_n'=\mathbb{C}[e_{2j+1}|0\leq j\leq \frac{n-1}{2}].
\non
\ena

\begin{prop}\label{prop-generation}
If $R_n$ is generated by 
$\{h_{i_1}\cdots h_{i_{n'}}|i_1,\cdots,i_{n'}\geq 0\}$ over $R_n'$,
it is generated by the same set over $R_n^{odd}$.
\end{prop}

To prove this proposition we notice the following lemma.

\begin{lemma}\label{sym-lem1}
In $R_n$ the following conditions are equivalent.

\noindent
$(1)$. $e_{2j+1}=0$ for $0\leq j\leq (n-1)/2$.

\noindent
$(2)$. $p_{2j+1}=0$ for $0\leq j\leq (n-1)/2$.
\end{lemma}
\vskip2mm

\noindent
{\it Proof}. 
By the formula ({\it cf}. \cite{Mac})
\bea
&&
ke_k=\sum_{r=1}^k(-1)^{r-1}p_re_{k-r},
\label{ep-rel}
\ena
each $e_k$ is expressed as a homogeneous polynomial of $p_1$,...,$p_k$
and vice versa where the degrees are defined by $\deg\, p_k=k$, 
$\deg\, e_k=k$.
Suppose that $e_j=0$ for all odd $j$. Then, for odd $k$,
if $p_k$ is written as a polynomial of $e_i$'s, each monomial contains
at least one $e_l$ with $l$ being odd and $p_k=0$. Thus (1)$\Rightarrow$(2)
is proved. The converse is similarly proved. \qed
\vskip5mm

\noindent
{\it Proof of Proposition \ref{prop-generation}}
\par
To avoid the notational complexity we give a proof in the case of $n$
even and use $2n$ instead of $n$. The case of $n$ odd is similarly
proved.

Let $E_{r_1\cdots r_n}=e_2^{r_1}e_4^{r_2}\cdots e_{2n}^{r_n}$.
Since $R_{2n}'=\mathbb{C}[e_1,e_3,\cdots,e_{2n-1}]$, 
$R_{2n}$ is a free $R_{2n}'$-module with the elements
$\{E_{r_1\cdots r_n}|r_j\in \mathbb{Z}_{\geq 0},\forall j\}$ as a basis.
By the assumption one can write
\bea
&&
E_{r_1\cdots r_n}=
\sum_{s_1,...,s_n\geq0}
f^{s_1\cdots s_n}_{r_1\cdots r_n}(e_1,e_3,\cdots,e_{2n-1})
h_{s_1}\cdots h_{s_n},
\non
\ena
for some $f^{s_1\cdots s_n}_{r_1\cdots r_n}\in R_{2n}'$.
Set
\bea
&&
F_{r_1\cdots r_n}=
\sum_{s_1,...,s_n\geq0}
f^{s_1\cdots s_n}_{r_1\cdots r_n}(p_1,p_3,\cdots,p_{2n-1})
h_{s_1}\cdots h_{s_n}.
\non
\ena

\begin{lemma}\label{sym-lem2}
The set of polynomials
$\{F_{r_1\cdots r_n}|r_j\in\mathbb{Z}_{\geq 0}\forall j\}$ is 
linearly independent over $R_{2n}^{odd}$.
\end{lemma}
\vskip2mm
\noindent
{\it Proof.}
It is sufficient to prove that $\{F_{r_1\cdots r_n}\}$ is linearly
independent over $\mathbb{C}$ at $p_1=p_3=\cdots=p_{2n-1}=0$.
By Lemma \ref{sym-lem1}
\bea
F_{r_1\cdots r_n}|_{p_{2j-1}=0, j\leq n}
&=&
\sum_{s_1,...,s_n\geq0}
f^{s_1\cdots s_n}_{r_1\cdots r_n}(0,\cdots,0)
(h_{s_1}\cdots h_{s_n})|_{p_{2j-1}=0, j\leq n}
\non
\\
&=&
E_{r_1\cdots r_n}|_{p_{2j-1}=0, j\leq n}.
\non
\ena
Thus we have to prove that 
$\tilde{E}_{r_1\cdots r_n}:=E_{r_1\cdots r_n}|_{p_{2j-1}=0, j\leq n}$
are linearly independent.

\begin{lemma}\label{sym-lem3}
At $e_{2j-1}=0, 1\leq j\leq n$ we have
\bea
&&
e_{2k}=-\frac{1}{2k}p_{2k}+\cdots,
\quad
1\leq k\leq n,
\non
\ena
where $\cdots$ part is a polynomial of $p_2$, $p_4,...,p_{2k-2}$.
\end{lemma}
\vskip2mm
\noindent
{\it Proof.}
The lemma easily follows from (\ref{ep-rel}). \qed

On the set $\{(r_1,\cdots,r_n)\}$ of non-negative integers we introduce
the lexicographical order from the right. Suppose that
\bea
&&
\sum c_{r_1\cdots r_n}\tilde{E}_{r_1\cdots r_n}=0,
\non
\ena
and $(r_1^0,\cdots,r_n^0)$ is the largest index such that 
$c_{r_1\cdots r_n}\neq 0$.
If one expresses each $\tilde{E}_{r_1\cdots r_n}$ as a polynomial
of $p_2,p_4,...,p_{2n}$, one has, by Lemma \ref{sym-lem3},
\bea
&&
c_{r_1^0\cdots r_n^0}p_2^{r_1^0}\cdots p_{2n}^{r_n^0}
+
\sum_{(r_1,\cdots,r_n)<(r_1^0,\cdots,r_n^0)}
c_{r_1 \cdots r_n}'p_2^{r_1}\cdots p_{2n}^{r_n}=0,
\non
\ena
for some constants $c_{r_1 \cdots r_n}'$.
Thus $c_{r_1^0 \cdots r_n^0}=0$ which contradicts the assumption.
Therefore Lemma \ref{sym-lem2} is proved. \qed

Since $R_{2n}=\oplus R_{2n}'E_{r_1\cdots r_n}$,
$\deg\, F_{r_1\cdots r_n}=\deg\, E_{r_1\cdots r_n}$ and 
$\ch\, R_{2n}'=\ch\, R_{2n}^{odd}$, we have
\bea
&&
\ch\,\Big( \oplus R_{2n}^{odd}F_{r_1\cdots r_n}\Big)=\ch\, R_{2n}.
\non
\ena
Thus $R_{2n}=\oplus R_{2n}^{odd}F_{r_1\cdots r_n}$ which proves
Proposition \ref{prop-generation}. \qed

\begin{prop}\label{prop-generation2}
The module $R_n$ is generated by 
$\{h_{i_1}\cdots h_{i_{n'}}|i_j\in\mathbb{Z}_{\geq 0}\forall j\}$
over $R_{n}'$.
\end{prop}
\vskip2mm
\noindent
{\it Proof.}
As in the previous proposition we give a proof for $n$ being even.
The proof of the case $n$ odd is similar. 
In this proof we again use $2n$ instead of $n$.
Then $(2n)'=n$.

For a partition $\lambda=(\lambda_1,\cdots,\lambda_m)$ we denote 
$s_\lambda$ the Schur function and $\ell(\lambda)$ the length
of $\lambda$ \cite{Mac}.

Let $L_{2n}$ be the submodule of $R_{2n}$ generated by 
$\{s_\lambda|\ell(\lambda)\leq n\}$ over $R_{2n}'$.
Due to the formula
\bea
&&
s_\lambda=\det(h_{\lambda_i-i+j})_{1\leq i,j\leq \ell(\lambda)},
\non
\ena
it is sufficient to prove $R_{2n}=L_{2n}$.
Since $\{s_\lambda|\ell(\lambda)\leq 2n\}$ is a $\mathbb{C}$-linear
basis of $R_{2n}$, we have to prove $s_\lambda\in L_{2n}$
for any partition $\lambda=(\lambda_1,\cdots,\lambda_{2n})$.
We prove this by induction on the degree 
$|\lambda|=\lambda_1+\cdots+\lambda_{2n}$ of $s_\lambda$.

If $|\lambda|\leq n$ for $\lambda=(\lambda_1,\cdots,\lambda_{2n})$,
$\lambda_j=0$ for $j\geq n+1$ and $s_\lambda\in L_{2n}$.

Let $d>n$. Suppose that $s_\mu\in L_{2n}$ for any partition 
$\mu=(\mu_1,\cdots,\mu_{2n})$ with $|\mu|<d$. 
We prove $s_\lambda\in L_{2n}$ for any 
$\lambda=(\lambda_1,\cdots,\lambda_{2n})$ with $|\lambda|=d$.
To this end we introduce the lexicographical order from the left 
on the set of partitions $\{(\mu_1,\cdots,\mu_{2n})\}$.
We prove $s_\lambda\in L_{2n}$ by descending induction on 
this order of $\lambda$.

For $\lambda=(d)$ we have $s_\lambda\in L_{2n}$ since $\ell(\lambda)=1\leq n$.

Let $\lambda$ be such that $\lambda<(d)$. Suppose that $s_\mu\in L_{2n}$
for any $\mu>\lambda$ with $|\mu|=d$. We assume that $\lambda$ is of the form
$\lambda=(\lambda_1,\cdots,\lambda_m)$, $\lambda_m\neq 0$.
Since $s_\lambda\in L_{2n}$ for $m\leq n$,
we assume that $n+1\leq m\leq 2n$.

First consider the case $\lambda_m=1$. Let 
$\tilde{\lambda}=(\lambda_1,\cdots,\lambda_{m-1})$. 
Then $|\tilde{\lambda}|=|\lambda|-1$ and $s_{\tilde{\lambda}}\in L_{2n}$
by the assumption of induction on $|\lambda|$.
Recall that $s_{(1^k)}=e_k$. Then $s_{\tilde{\lambda}}s_{(1)}\in L_{2n}$.
By Pieri's formula ({\it cf.} \cite{Mac})
\bea
&&
s_{\tilde{\lambda}}s_{(1)}=s_\lambda+\sum_{\mu>\lambda}c_\mu s_\mu,
\label{decomp1}
\ena
for some constants $c_\mu$.
By the hypothesis of induction on the order of $\lambda$ the second term
of the right hand side of (\ref{decomp1}) is in $L_{2n}$.
Thus $s_\lambda\in L_{2n}$.

Next consider the case $m=2n$.
We set $\tilde{\lambda}=(\lambda_1,\cdots,\lambda_{2n-1},\lambda_{2n}-1)$.
Then $\tilde{\lambda}$ is a partition and $s_{\tilde{\lambda}}\in L_{2n}$
by the induction hypothesis on $|\lambda|$. 
Then $s_{\tilde{\lambda}}s_{(1)}\in L_{2n}$ and the equation (\ref{decomp1})
holds. Thus $s_\lambda\in L_{2n}$.

Let us assume $\lambda_m\geq 2$ and $m<2n$.

\begin{lemma}\label{sym-lem4}
Let $m\leq m'\leq 2n$ and $i_j\in\mathbb{Z}_{\geq0}$, $m+m'-2n+1\leq j\leq m$ be such that
$\sum\, i_j=2n-m'$ and
\bea
&&
\lambda^{(m')}:=
(\lambda_1,\cdots,\lambda_{m+m'-2n},\lambda_{m+m'-2n+1}-i_{m+m'-2n+1},
\cdots,\lambda_m-i_m,1^{2n-m'})
\non
\ena
is a partition. Then $s_{\lambda^{(m')}}\in L_{2n}$.
\end{lemma}
\vskip2mm
\noindent
{\it Proof.}
In this proof we use the following notations.
Let $\ep_i=(0,\cdots,1,\cdots,0)$ where $1$ is on the $i$-th position.
For $N\geq 1$ and $a=(a_i),b=(b_i)\in \mathbb{Z}^N$ we define $a+b=(a_i+b_i)$.
For $N_1\leq N_2$, $a=(a_i)\in \mathbb{Z}^{N_1}$ is considered as an element of
$\mathbb{Z}^{N_2}$ by setting $a=(a_i)$, $a_i=0$ for $N_1<i$.

We prove the lemma by induction on $m'$.

Suppose that $m'=m$. Set
\bea
&&
\tilde{\lambda}^{(m)}:=
(\lambda_1,\cdots,\lambda_{2m-2n-1},\lambda_{2m-2n}-1,
\lambda_{2m-2n+1}-i_{2m-2n+1}-1,
\cdots,\lambda_m-i_m-1).
\non
\ena
Since $m<2n$ and $\lambda^{(m)}$ is a partition, $\lambda_m-i_m\geq 1$.
Then $\tilde{\lambda}^{(m)}$ becomes a partition and 
$|\tilde{\lambda}^{(m)}|=|\lambda|-(2(2n-m)+1)<|\lambda|$.
Thus $s_{\tilde{\lambda}^{(m)}}\in L_{2n}$ by the assumption of induction
on $|\lambda|$. Then
\bea
&&
s_{\tilde{\lambda}^{(m)}}s_{(1^{2(2n-m)+1})}
=s_{\lambda^{(m)}}+\sum_{\mu>\lambda}c_\mu s_\mu,
\non
\ena
for some constants $c_\mu$'s.
Thus $s_{\lambda^{(m)}}\in L_{2n}$ by induction on the order of $\lambda$.

Assume that $m<m'\leq 2n$ and that the lemma holds for any $m''$
satisfying $m\leq m''<m'$.
Let 
\bea
&&
\tilde{\lambda}^{(m')}
=
(\lambda_1,\cdots,\lambda_{m+m'-2n-1},\lambda_{m+m'-2n}-1,
\lambda_{m+m'-2n+1}-i_{m+m'-2n+1}-1,
\cdots,\lambda_m-i_m-1).
\non
\ena
Here $\lambda_m-i_m\geq 1$. In fact if $m'<2n$ then it is obvious that $\lambda_m-i_m\geq 1$. If $m'=2n$ then $i_j=0$ for all $j$ and 
$\lambda_m-i_m=\lambda_m\geq 2$ by assumption.
Then $\tilde{\lambda}^{(m')}$ is a partition.
Since $|\tilde{\lambda}^{(m')}|<|\lambda|$, 
$s_{\tilde{\lambda}^{(m')}}\in L_{2n}$
by the hypothesis of induction on $|\lambda|$.
Then
\bea
&&
s_{\tilde{\lambda}^{(m')}}s_{(1^{2(2n-m')+1})}
=s_{\lambda^{(m')}}
+\sum_{k\geq 1}\sum_{\{j_l\}}s_{\lambda^{[k]}(\{j_l\})}
+\sum_{\mu>\lambda}c_\mu s_\mu.
\non
\ena
Here $\{j_l\}$ and $\lambda^{[k]}(\{j_l\})$ are as follows.
For a given $k$, $\{j_l\}$ is a set of numbers satisfying
\bea
&&
m+m'-2n\leq j_1<\cdots<j_{2n-m'-k+1}\leq m.
\non
\ena
Then $\lambda^{[k]}(\{j_l\})$ is defined by
\bea
&&
\lambda^{[k]}(\{j_l\})=
\tilde{\lambda}^{(m')}+\sum_{l}\ep_{j_l}+(0^{m},1^{2n-m'+k}).
\non
\ena
If $\lambda^{[k]}(\{j_l\})$ is not a partition we define 
$s_{\lambda^{[k]}(\{j_l\})}=0$.
If $\lambda^{[k]}(\{j_l\})$ is a partition and 
$\ell\big(\lambda^{[k]}(\{j_l\})\big)>2n$ then 
$s_{\lambda^{[k]}(\{j_l\})}=0$.
Consider $k$ such that $\lambda^{[k]}(\{j_l\})$ is a partition and 
$\ell\big(\lambda^{[k]}(\{j_l\})\big)\leq 2n$.
Since
\bea
&&
\lambda=\tilde{\lambda}^{(m')}+\sum_{p=m+m'-2n}^m \ep_p
+\sum_{j=m+m'-2n+1}^m i_j\ep_j,
\non
\ena
we have
\bea
&&
\lambda^{[k]}(\{j_l\})-\lambda
=-\sum_{p=m+m'-2n}^m \ep_p
+\sum_{l=1}^{2n-m'-k+1} \ep_{j_l}
-\sum_{j=m+m'-2n+1}^m i_j\ep_j
+(0^{m},1^{2n-m'+k}).
\label{lambda-diff}
\ena
Let us write
\bea
&&
\sum_{p=m+m'-2n}^m \ep_p
-\sum_{l=1}^{2n-m'-k+1} \ep_{j_l}
+\sum_{j=m+m'-2n+1}^m i_j\ep_j
=
(0^{m+m'-2n-k},i'_{m+m'-k-2n+1},\cdots,i'_m)=:\mathbf{i}'.
\non
\ena
Then $i'_j\geq 0$ for any $j$, $\sum\, i_j'=2n-(m'-k)$ and 
\bea
&&
\lambda-\mathbf{i}'+(0^{m},1^{2n-m'+k})=\lambda^{[k]}(\{j_l\})
\non
\ena
is a partition.
Thus $(i_j')$ satisfies the condition of $(i_j)$ for $m'-k$.
 Since $k\geq 1$, 
$s_{\lambda^{[k]}}(\{j_l\})\in L_{2n}$ by the assumption of induction on $m'$.
Thus $s_{\lambda^{(m')}}\in L_{2n}$. \qed
\vskip2mm

Consider the case $m'=2n$ in Lemma \ref{sym-lem4}.
Then $\lambda^{(2m)}=\lambda$. Thus $s_\lambda\in L_{2n}$.
This completes the proof of Proposition \ref{prop-generation2}.
\qed

\begin{lemma}\label{sym-lem5}
If $e_{2j+1}=0$ for all $j\leq (n-1)/2$, $h_{2r-1}=0$ for any $r\geq 1$.
\end{lemma}
\vskip2mm
\noindent
{\it Proof.}
The lemma is easily proved by induction on $r$ using the relation
\bea
&&
\sum_{i=0}^m(-1)^{i}e_ih_{m-i}=0,
\quad
m\geq 1.
\non
\ena
\qed
\vskip5mm

\noindent
{\it Proof of Theorem \ref{th-symfunc}}
\par
Consider the specialization $p_{2j+1}=0$ for all $j\leq (n-1)/2$.
Then 
$R_n=\mathbb{C}[p_2,p_4,\cdots,p_{2n'}]$ and 
\bea
&&
\ch\, R_n=\prod_{j=1}^{n'}\frac{1}{1-q^{2j}}.
\label{chRn-1}
\ena
By Lemma \ref{sym-lem1} we have $e_{2j+1}=0$ for $j\leq (n-1)/2$.
Then $h_{2r-1}=0$ for any $r\geq 1$ by Lemma \ref{sym-lem5}.
Thus
\bea
&&
R_n=\sum_{0\leq r_1\leq\cdots\leq r_{n'}}\mathbb{C}
h_{2r_1}\cdots h_{2r_{n'}},
\label{chRn-2}
\ena
by Proposition \ref{prop-generation} and Proposition \ref{prop-generation2}.
For two formal power series $f(q)$, $g(q)$ of $q$ with the coefficients
in $\mathbb{R}$ we denote $g(q)\leq f(q)$ if all coefficients
of $f(q)-g(q)$ is non-negative. Then
\bea
&&
\ch(\text{RHS of (\ref{chRn-2})})
\leq
\sum_{0\leq r_1\leq \cdots\leq r_{n'}}q^{2(r_1+\cdots+r_{n'})}
=\prod_{j=1}^{n'}\frac{1}{1-q^{2j}}=\ch\, R_n.
\label{inequality}
\ena
Thus the inequality in (\ref{inequality}) is in fact the equality.
This means that 
$\{h_{2r_1}\cdots h_{2r_{n'}}|0\leq r_1\leq \cdots\leq r_{n'}\}$
are linearly independent over $\mathbb{C}$ at $p_{2j+1}=0$, $j\leq (n-1)/2$.
Thus $\{h_{2r_1}\cdots h_{2r_{n'}}|0\leq r_1\leq \cdots\leq r_{n'}\}$
are linearly independent over $R_n^{odd}$. 
Then
\bea
&&
\ch\Big( \oplus _{0\leq r_1\leq \cdots\leq r_{n'}}
R_n^{odd}h_{2r_1}\cdots h_{2r_{n'}}\Big)
=\prod_{j=1}^n\frac{1}{1-q^j}=\ch\, R_n.
\non
\ena
This completes the proof.
\qed

\section{Anti-chiral subspace}
In this section we construct local operators whose values of spins are minus
of those of the operators in the chiral subspace in section 3.
We call the subspace of such operators the anti-chiral subspace.

Let $R_n^{-}=\mathbb{C}[x_1^{-1},\cdots,x_n^{-1}]^{S_n}$, 
$H^{(n)-}=\oplus_{j=0}^{n-1}R_n^{-}X^{-j}$, 
$\prod_{j=1}^n(1+x_j^{-1}t)=\sum_{j=0}^ne^{(n)-}_j t^j$ and
\bea
&&
U_{n,\ell}^{-}=\{P\in \wedge^\ell H^{(n)-}|\rho_{\pm}(P)=0\}.
\non
\ena
Define the map
\bea
&&
{}^{-}:\wedge^\ell H^{(n)\pm}\lar \wedge^\ell H^{(n)\mp},
\non
\ena
by
\bea
&&
P(X_1,\cdots,X_\ell|x_1,\cdots,x_n)
\mapsto
P^{-}:=P(X_1^{-1},\cdots,X_\ell^{-1}|x_1^{-1},\cdots,x_n^{-1}),
\non
\ena
where we set $H^{(n)+}=H^{(n)}$. Obviously $P^{--}=P$ and the map ${}^{-}$ 
gives an isomorphism of $\mathbb{C}$-vector spaces.
This map induces an isomorphism between $U_{n,\ell}$ and $U_{n,\ell}^{-}$
as a vector space. In fact for $P\in \wedge^\ell H^{(n)}$ we have
\bea
&&
P|_{X_\ell=\pm x^{-1}, x_n=-x_{n-1}}=0
\Leftrightarrow
P^{-}|_{X_\ell^{-1}=\pm x, x_n^{-1}=-x_{n-1}^{-1}}=0.
\non
\ena
The last condition is equivalent to $\rho_\pm(P^{-})=0$.
Moreover the map satisfies, for $f\in R_n$ and $P\in \wedge^\ell H^{(n)}$,
\bea
&&
(fP)^{-}=f^{-}P^{-}.
\non
\ena
Thus we have

\begin{cor}
The module $U_{n,\ell}^{-}$ is a free $R_n^{-}$-module with 
the elements $\{v^{(2n)-}_I\wedge w^{(2n)-}_J\wedge \xi^{(2n)-}_K\}$
as a basis, where $I,J,K$ satisfy the same conditions as in Theorem \ref{basis}.\end{cor}

We set
\bea
&&
M_{n,\ell}^{-}=
\frac{U_{n,\ell}^{-}}
{U_{n,\ell-1}^{-} \wedge \Xi^{(n)-}_1
+U_{n,\ell-2}^{-}\wedge \Xi^{(n)-}_2},
\non
\ena
and, for $i=0,1$, $\lambda\in \mathbb{Z}_{\geq 0}$,
\bea
&&
M^{(i)-}_{\lambda}=
\oplus_{n\equiv i\, mod.\, 2,\, n-2\ell=\lambda} M_{n,\ell}^{-}.
\non
\ena
We introduce the degree on $M_{n,\ell}^{-}$ and $M^{(i)-}_{\lambda}$
by $\deg_3\, P=-n^2/4+\deg_1\, P$.
It is obvious that $\deg_1\, P^{-}=-\deg_1\, P$ for a homogeneous element 
$P\in H^{(n)\pm}$.
Then
\bea
&&
\ch\, M^{(i)-}_{\lambda}=
\sum_{n-2\ell=\lambda,\, n\equiv i\, mod.\, 2}
\frac{q^{-\frac{n^2}{4}}}{[n]_{q^{-1}}!}
\Big(
\qibc{n}{\ell}-\qibc{n}{\ell-1}
\Big).
\non
\ena

In the remaining part of this section we consider the case of 
$n$ even and use $2n$ instead of $n$.
Let us consider the part of $\Psi_P$ which is relevant to determining
the homogeneous degree of $\Psi_P$.
It is
\bea
&&
\exp\Big(\frac{n}{2}\sum_{j=1}^{2n}\beta_j\Big)
\frac{P}{\prod_{a=1}^\ell\prod_{j=1}^{2n}(1-X_a x_j)}.
\non
\ena
We rewrite it as
\bea
&&
\exp\Big(-\frac{n}{2}\sum_{j=1}^{2n}\beta_j\Big)
\frac{(e^{(2n)}_{2n})^{n-\ell}(\prod_{a=1}^\ell X_a^{-2n}) P}
{\prod_{a=1}^\ell\prod_{j=1}^{2n}(1-X_a^{-1} x_j^{-1})}.
\non
\ena
For a polynomial $P=P(X_1,\cdots,X_\ell|x_1,\cdots,x_{2n})$ of
$X_a$'s such that $\deg_{X_a}P\leq 2n$ set
\bea
&&
\widehat{P}=(e^{(2n)}_{2n})^{n-\ell}(\prod_{a=1}^\ell X_a^{-2n}) P.
\non
\ena
It is a polynomial of $X_a^{-1}$'s such that $\deg_{X_a^{-1}}\, P\leq 2n$.
By this notation it is obvious that $\Psi_P$ is homogeneous of degree
$-n^2+\deg_1\widehat{P}=\deg_3 \widehat{P}$.
Rewriting Theorem \ref{nt-th} in terms of $\widehat{P}$ we have

\begin{theorem}\label{anti-th1}
Let $\widehat{P}_{2n}(X_1,...,X_{\ell_{2n}}|x_1,...,x_{2n})$, $n\geq m$ 
be polynomials 
of $X_a^{-1}$ with the coefficients in the ring of symmetric 
Laurant polynomials of $x_j$'s such that 
$\deg_{X_a^{-1}}\, \widehat{P}_{2n}\leq 2n$
for $1\leq a\leq \ell_{2n}$. 
Suppose that $(\widehat{P}_{2n})_{n=0}^\infty$ satisfy 
the following conditions: there exists a set of polynomials 
$(\widehat{P}_{2n}')$ 
such that
\bea
&&
Asym\Big(\overline{\widehat{P}_{2n}}\Big)=
Asym\Big(\prod_{a=1}^{\ell_{2n}-1}(1-x^{-2}X_a^{-2})\widehat{P}_{2n}'\Big),
\label{anti-cond1}
\\
&&
\widehat{P}_{2n}'|_{X_{\ell_{2n}}=\pm x^{-1}}
=\pm(-1)^{n-1} x^{2n-1}d_{2n}\widehat{P}_{2n-2},
\label{anti-cond2}
\ena
where the notations are same as those in Theorem \ref{nt-th}.
Set $P_{2n}=(e^{(2n)}_{2n})^{r-m}(\prod_{a=1}^{\ell_{2n}} X_a^{2n})
\widehat{P}_{2n}$. 
Then $(\Psi_{P_{2n}})_{n=0}^\infty$
satisfies $(\ref{iff-1})$, $(\ref{iff-2})$, $(\ref{iff-3})$.
\end{theorem}
\vskip3mm

\begin{prop}\label{anti-ff}
The following $(\widehat{P}_{2n})_{n=0}^\infty$ satisfies the 
conditions of Theorem \ref{anti-th1}:
\bea
&&
\widehat{P}_{2m}=v^{(2m)-}_I\wedge w^{(2m)-}_J\wedge \xi^{(2m)-}_K,
\non
\\
&&
\widehat{P}_{2n}=
c_{2n}^{-}
v^{(2n)-}_I\wedge w^{(2n)-}_J\wedge \xi^{(2n)-}_K
\prod_{a=r+1}^{\ell_{2n}}
X_a^{-(2n+1+2r-2a)},
\non
\\
&&
c_{2n}^{-}=
(i\zeta(-\pi i))^{m-n}
(2\pi i)^{(m-n)(n+r)},
\non
\ena
where $\ell_1+\ell_2+2\ell_3=r$ and $P_{2n}=0$, $n<m$.
\end{prop}

It is possible to multiply $Q^{(-)}_{2n}(z)$ to $\widehat{P}_{2n}$ 
simultaneously
without destroying the conditions in Theorem \ref{anti-th1} if the degrees
in $X_a^{-1}$ remains the permitted range.
Let us write the final formula of $P$ corresponding to the operators
in the anti-chiral subspace.
Define
\bea
&&
P_{2n}^{(-)}=
c_{2n}^{-}
(e_{2n}^{(2n)})^{-m+r}
E_{odd}^{(2n)}(t)
\prod_{i=1}^m Q^{(-)}_{2n}(z_i)
v^{(2n)-}_I\wedge w^{(2n)-}_J\wedge \xi^{(2n)-}_K
\non
\\
&&
\qquad\qquad
\times
\prod_{a=1}^r X_a^{2n}
\prod_{a=r+1}^{\ell_{2n}}
X_a^{2(a-r)-1},
\label{anti-P}
\ena
and, as in the chiral case,
\bea
&&
P_{2n}^{(-)}=\sum P_{2n,\alpha,\gamma}^{(-)}t^\alpha z^\gamma.
\non
\ena
Then
\bea
&&
\deg_3\, P_{2n,\alpha,\gamma}^{(-)}
=-m^2+\sum i\,\alpha_i-\sum \gamma_i-
\deg_1(v^{(2m)}_I\wedge w^{(2m)}_J\wedge \xi^{(2m)}_K),
\non
\ena
for all $n\geq m$.

\begin{cor}
$(1)$. The set of functions $(\Psi_{P_{2n}^{(-)}})_{n=0}^\infty$ satisfies
$(\ref{iff-1})$, $(\ref{iff-2})$, $(\ref{iff-3})$ and each 
$\Psi_{P_{2n}^{(-)}}$ takes the value in 
$(V^{\ot 2n})^{sing}_{\lambda}$, $\lambda=2m-2r$.

\noindent
$(2)$. The following property holds for all $n$:
\bea
&&
\Psi_{P_{2n,\alpha,\gamma}^{(-)}}(\beta_1+\theta,\cdots,\beta_{2n}+\theta)
=\exp(\theta\deg_3\, P_{2m,\alpha,\gamma}^{(-)})
\Psi_{P_{2n,\alpha,\gamma}^{(-)}}(\beta_1,\cdots,\beta_{2n}).
\non
\ena
\end{cor}

Obviously the following equation holds:
\bea
&&
\sum_{\alpha,\gamma}\mathbb{C}P_{2m,\alpha,\gamma}^{(-)}=U_{2m,r}^{(-)},
\non
\ena
where the summation is taken for all $\gamma=(\gamma_1,\gamma_2,\cdots)$ and
$\alpha=(\cdots,\alpha_{-3},\alpha_{-1},\alpha_1,\cdots)$ with 
$\alpha_i=0$, $i\geq 1$.

\section{Concluding remarks}
In \cite{Smir3} Smirnov has given a systematic construction of form
factors of a large number of chargeless local operators in Sine-Gordon
model. This construction is extended to the charged operators in SU(2) ITM 
in \cite{NT}.
In these descriptions of local operators, however, finding the subspace 
which has the same character, with respect to spins, as the chiral 
subspace of local operators in the corresponding conformal field theory 
has not been successful up to now.

In this paper we have given an alternative construction of local operators 
in SU(2) ITM based on the results of our previous paper \cite{N} 
and the ideas of \cite{Smir3,NT}. 
With this description of local operators we have
specified the subspace of operators isomorphic to the level one integrable
highest weight representation of $\widehat{sl_2}$ which is the chiral space
of the level one $su(2)$ WZW model.

Let us give a brief comment here on the differences of the formulae 
(\ref{ff-1}), (\ref{ff-2}) and those of Smirnov in \cite{Smir3}.

In \cite{Smir3} the initial condition is taken as
\bea
&&
P_{2m}=(e^{(2m)}_{2m})^{-N}f\Delta^{+}_{2m}\prod_{a=1}^rX_a^{r_a},
\quad
0\leq r_a\leq 2m-1,
\quad
\Delta^{+}_{2m}=\prod_{i<j}(x_i+x_j),
\non
\ena
for some $N\geq 0$ and $f\in R_{2m}$.
Since $\overline{\Delta^{+}_{2m}}=0$, 
$Asym_{X_1,\cdots,X_r}\Big((e^{(2m)}_{2m})^{N}P_{2m}\Big)\in U_{2m,\ell}$.
The important point is that, in general, $U_{2m,r}$ is not generated 
by such functions only.
That is why we have introduced the polynomials $v^{(n)}_i$, $w^{(n)}_j$
and $\xi^{(n)}_k$. This difference is crucial in the calculation of the 
character of the chiral subspace of local operators.

Due to this difference of the initial condition the degrees of freedom corresponding to the parameter $t_{2s}$, $s\in \mathbb{Z}$ in \cite{Smir3} could not be
introduced in our formula. This drawback is compensated by the introduction
of the functions $Q_{2n}^{(\pm)}(z)$. 
With this respect the structure of symmetric polynomials given in 
Theorem \ref{th-symfunc} is important.

\appendix
\section{Function $\zeta(\beta)$}
The function $\zeta(\beta)$ was introduced in \cite{KS,Smir1,Smir2}.
It is conveniently described using the double gamma function of 
Barnes \cite{NT,JM}. We review it here.
We follow \cite{Sh} about the notations of double gamma function.

For a set of positive real numbers $\omega=(\omega_1,\omega_2)$
define
\bea
&&
\Gamma_2(z|\omega)^{-1}=
z\exp\Big(
\gamma_{22}(\omega)+\frac{z^2}{2}\gamma_{21}(\omega)
\Big)
\hbox{$\prod'$}
\Big(
1+\frac{z}{m\omega_1+n\omega_2}
\Big)
\exp\Big(
-\frac{z}{m\omega_1+n\omega_2}+\frac{z^2}{2(m\omega_1+n\omega_2)^2}
\Big)
,
\non
\ena
where the product $\prod'$ is over all sets of non-negative integers
except $(m,n)=(0,0)$ and $\gamma_{22}(\omega)$, $\gamma_{22}(\omega)$
are given by the formulae
\bea
&&
-\gamma_{22}(\omega)=
\frac{1}{\omega_1}\sum_{n=1}^\infty
\Big(
\psi(\frac{n\omega_2}{\omega_1})-\log(\frac{n\omega_2}{\omega_1})
+\frac{\omega_1}{2n\omega_2}
\Big)
+\frac{1}{2}\big(\frac{1}{\omega_1}+\frac{1}{\omega_2}\big)\log\omega_1
\non
\\
&&
-\frac{1}{2\omega_1}(\gamma-\log\, 2\pi)
+\frac{\omega_1-\omega_2}{2\omega_1\omega_2}
\log\frac{\omega_2}{\omega_1}
-\frac{\gamma}{2}\big(\frac{1}{\omega_1}+\frac{1}{\omega_2}\big),
\label{gamma22}
\\
&&
-\gamma_{21}=
\frac{1}{\omega_1^2}
\Big(
\psi'(\frac{n\omega_2}{\omega_1})
-\frac{\omega_1}{n\omega_2}
\Big)
+\frac{\pi^2}{6\omega_1^2}
-\frac{1}{\omega_1\omega_2}\log\omega_2
+\frac{\gamma}{\omega_1\omega_2}.
\label{gamma21}
\ena
Here $\psi$ is the logarithmic derivative of the gamma function and
$\gamma$ is the Euler's constant.
Then $\Gamma_2(z|\omega)^{-1}$ becomes an entire function of $z$,
is symmetric with respect to $\omega_1$ and $\omega_2$ and satisfies
the following difference equation
\bea
&&
\frac{\Gamma_2(z+\omega_1|\omega)}{\Gamma_2(z|\omega)}
=
\sqrt{2\pi}
\exp\Big((-\frac{z}{\omega_2}+\frac{1}{2})\log\omega_2\Big)
\frac{1}{\Gamma(\frac{z}{\omega_2})}.
\label{diffeq-gamma2}
\ena
Set $\Gamma_2(z)=\Gamma_2(z|2\pi,2\pi)$ and define
\bea
&&
\zeta(\beta)=
\frac{\Gamma_2(-i\beta+\pi)\Gamma_2(i\beta+3\pi)}
{\Gamma_2(-i\beta)\Gamma_2(i\beta+2\pi)}.
\label{zeta}
\ena
Then it satisfies the following equations
\bea
&&
\zeta(\beta-2\pi i)=\zeta(-\beta),
\label{zeta-eq1}
\\
&&
\zeta(\beta)\zeta(\beta-\pi i)
=
\frac{(2\pi )^{3/2}}
{\Gamma(\frac{-i\beta+\pi}{2\pi})\Gamma(\frac{i\beta}{2\pi})},
\label{zeta-eq2}
\\
&&
\frac{\zeta(-\beta)}{\zeta(\beta)}
=\frac{\Gamma(\frac{\pi i+\beta}{2\pi i})\Gamma(\frac{-\beta}{2\pi i})}
{\Gamma(\frac{\pi i-\beta}{2\pi i})\Gamma(\frac{\beta}{2\pi i})}.
\label{zeta-eq3}
\ena
The right hand side of (\ref{zeta-eq3}) is $S_0(\beta)$ in the $S$-matrix.

\end{document}